\documentstyle{ptptex}

\input amssymb.sty
\newtheorem{thm}{Theorem}[section]

\newtheorem{ex}[thm]{Example}

\newcommand{\reals}{{\Bbb R}}

\newcommand{\Ad}{{\rm Ad}}

\newcommand{\cala}{{\cal A}}
\newcommand{\cald}{{\cal D}}

\newcommand{\call}{{\cal L}}

\newcommand{\calu}{{\cal U}}

\newcommand{\del}{\partial}

\newcommand{\half}{\textstyle{\frac{1}{2}}}

\newcommand{\frakg}{\mathfrak{g}}

\newcommand{\arrows}{\,\lower1pt\hbox{$\longrightarrow$}\hskip-.24in\raise2pt
             \hbox{$\longrightarrow$}\,}

\newcommand{\Severa}{\v{S}evera~}
\newcommand{\Klimcik}{Klim\v c\'\i k}
\newcommand{\wt}{\widetilde}
\title{{\bf Poisson geometry\\with a 3-form background}}
\author{Pavol \v Severa\thanks{Research supported by the European Postdoctoral Institute (EPDI).}\\
Department of Theoretical Physics\\Comenius University\\84215 Bratislava, Slovakia
\\{\small(severa@sophia.dtp.fmph.uniba.sk)}
\\and Alan
Weinstein\thanks{Research partially supported by NSF Grant
DMS-99-71505. \newline{\it MSC2000 Subject
Classification:} Primary 53D17, Secondary 58H05, 58H15
\newline{\it Keywords}: Poisson structure, Dirac structure, Lie
algebroid, Courant algebroid, deformation quantization}
\\Department of Mathematics\\ University of California\\ Berkeley, CA
94720 USA\\ {\small(alanw@math.berkeley.edu)}}
\abst{We study a modification of Poisson geometry by a closed 3-form.
Just as for ordinary Poisson structures, these ``twisted" Poisson
structures are conveniently described as Dirac structures in
suitable Courant algebroids. The additive group of 2-forms acts on
twisted Poisson structures and permits them to be seen as glued
from ordinary Poisson structures defined on local patches. We
conclude with remarks on deformation quantization and twisted
symplectic groupoids.}

\begin{document}
\maketitle
The ideas presented in this note grew out of an attempt to understand
how Poisson geometry on  a manifold is affected by the presence of
a closed 3-form ``field''.  Such forms  are playing an important
role in contemporary string theory.  We refer, for example, to
Park \cite{pa:topological} as well as to
Cornalba and Sciappa \cite{co-sc:nonassociative} 
and \Klimcik ~and Str\"obl \cite{kl-st:wzw}. 
Our aim here is to show that the notions of Courant
algebroid and Dirac structure provide a framework in which one can
easily carry out computations in Poisson geometry in the presence
of a background 3-form. It seems clear that a proper understanding
of the global effect of such a 3-form involves gerbes (see for
example Brylinski \cite{br:loop}); our work here should at least
partially substantiate the claim that Courant algebroids are
appropriate infinitesimal objects to associate with gerbes.

Our work was stimulated in part by the many talks at the Workshop
on Deformation Quantization and String Theory at Keio University
(March, 2001) in which such 3-forms played an essential role.  It is
essentially an application of some of the ideas contained in
a series of letters from \Severa to Weinstein written in
1998.  Some of the material in this paper was presented in June,
2001 at the
Colloque en l'honneur d'Yvette Kosmann-Schwarzbach at the Institut
Henri Poincar\'e and the Workshop on Poisson Geometry at
the Erwin Schr\"odinger Institute.

Our basic idea is as follows.  Poisson structures on a manifold $M$
may be identified with certain Dirac structures in the standard
Courant algebroid $E_0=TM\oplus T^*M$; see Courant \cite{co:dirac}
and Liu {\it et al.}~\cite{li-we-xu:manin} It turns out that a
closed 3-form $\phi$ on $M$ may be used to modify the bracket on
$E_0$, yielding a new Courant algebroid $E_{\phi}$.  A bivector
$\pi$ on $M$ now corresponds to a Dirac structure in $E_\phi$ if
and only if it satisfies the equation
\begin{equation}
\label{eq-quasi}
[\pi,\pi]=\wedge^3\tilde{\pi}(\phi).
\end{equation}
Here, and elsewhere in this note,
the operator $\tilde{B}:V\to V^*$ is   defined for any bilinear
form $B$ on a vector space $B$ by
$\tilde{B}(\alpha)(\beta)=B(\alpha,\beta)$ for $\alpha$ and $\beta$ in $V$.
Given $\phi$, we refer to  solutions $\pi$
of (\ref{eq-quasi}) as {\bf $\phi$-Poisson
structures} or, if we do not want to specify $\phi$, {\bf twisted Poisson
structures}.

Twisted Poisson structures arose from the study of topological sigma
 models in the work of Park
\cite{pa:topological}, as well as 
\Klimcik{} and Str\"obl
 \cite{kl-st:wzw}, who called them {\bf WZW-Poisson
structures}.  A related notion of {\bf quasi-Poisson
structure} on a manifold with a group action has been introduced
by Alekseev {\em et al.}
\cite{al-ko-me:quasi}, in connection with the theory of group-valued
momentum mappings.  Their Jacobi anomaly comes from a trivector on
the group rather than a 3-form on the manifold, but it is possible
that there is a general notion (perhaps involving Lie algebroid
actions) which will incorporate both twisted Poisson and quasi-Poisson
 structures. 

We would like to thank Anton Alekseev, Ctirad \Klimcik, Yvette
Kosmann-Schwarzbach, Jiang-hua Lu, Dmitry Roytenberg, Thomas Str\"obl,
and Ping Xu for helpful advice.  We thank Jae-Suk Park for calling his
earlier work on this subject to our attention.  Finally,
A.W. would also like to thank
colleagues at Keio University, \'Ecole Polytechnique, Institut de
Math{\'e}matiques de Jussieu, and Erwin Schr\"odinger Institute
for their hospitality during various stages of this work.

\section{Courant algebroids and Dirac structures}

The original bracket of Courant \cite{co:dirac} and
the Courant algebroid brackets of Liu {\it et al.}~\cite{li-we-xu:manin} were
skew-symmetric.  A non-skew-symmetric version of the bracket was also
introduced by Liu {\it et al.}~\cite{li-we-xu:manin}, where some of its nice properties
were noted.  It was then observed by Kosmann-Schwarzbach, Xu, and the
first author (all unpublished) that this bracket satisfied the Jacobi
identity written in Leibniz form.  We will use the non-skew-symmetric
bracket in this paper.

A {\bf Courant algebroid} over a manifold $M$
is a vector bundle $E\to M$
equipped with a field of nondegenerate  symmetric  bilinear forms
 $( \cdot , \cdot )$ on the fibres,  an
$\reals$-bilinear bracket
$[\cdot , \cdot ]: \Gamma(E)\times\Gamma(E)\rightarrow\Gamma(E)$
 on the space of sections of $E$,
and a bundle map $\rho :E\to TM$ (the {\bf anchor}),
 such that the following
properties are  satisfied:
 \begin{enumerate}
\item for any $e_{1}, e_{2}, e_{3}\in \Gamma (E)$,
$[e_{1},[e_{2}, e_{3}]]=[[e_1,e_2],e_3]+[e_2,[e_1,e_3]];$
\item  for any $e_{1}, e_{2} \in \Gamma (E)$,
$\rho [e_{1}, e_{2}]=[\rho e_{1}, \rho  e_{2}];$
\item  for any $e_{1}, e_{2} \in \Gamma (E)$ and $f\in C^{\infty} (M)$,
$[e_{1}, fe_{2}]=f[e_{1}, e_{2}]+(\rho (e_{1})f)e_{2} ;$
\item for any $e, h_{1}, h_{2} \in \Gamma (E)$,
  $\rho (e) (h_{1}, h_{2})=([e , h_{1}] ,
h_{2})+(h_{1}, [e , h_{2}] )$;
\item for any $e\in\Gamma(E)$, $[e,e]=\cald(e,e)$,
\end{enumerate}
where
$\cald :  C^{\infty}(M)\to \Gamma (E)$
is  the map   defined
  by $\cald = \half \beta^{-1}\rho^{*} d$,
 where $\beta $ is
the isomorphism between $E$
and $E^*$ given by the bilinear form. In other words,
\begin{equation}
\label{eq:D}
(\cald f , e)= \half  \rho (e) f .
\end{equation}

Equivalently, instead of the bracket $[\cdot,\cdot]$, we can use a
linear map $e\mapsto Z_e$ which maps sections of $E$ to vector
fields on the total space of $E$. The vector field $Z_e$ is a lift
of $\rho(e)$ from $M$ to $E$, and the first four axioms just say
that the flows of $Z_e$'s preserve the structure of $E$. The
bracket $[e_1,e_2]$ is the Lie derivative of $e_2$ by $Z_{e_1}$.
This feature of the present non-skew-symmetric bracket makes it
particularly convenient.

A {\bf Dirac structure} in $E$, also called an $E$-Dirac structure
on $M$, is a maximal isotropic subbundle $L$ of $E$ whose sections are
closed under the bracket, i.e. which is  preserved by
the flow of $Z_e$ for any $e\in \Gamma(L)$. The properties of a
Courant algebroid imply that the restriction of the  bracket and
anchor to any Dirac structure $L$ form a Lie algebroid structure on
$L$.

On any manifold $M$, we have the standard Courant algebroid
$E_0=TM\oplus T^*M$ with bilinear form
$((X_1,\xi_1),(X_2,\xi_2))=\xi_1(X_2)+\xi_2(X_1)$, anchor $\rho(X,\xi)=X$, and
 bracket
\begin{equation}
\label{eq-courant}
[(X_1, \xi_1),( X_2,\xi_2)]=
([X_1,X_2], {\cal L}_{X_1}\xi_2-i_{X_2}d\xi_1).
\end{equation}
Thus, $Z_{(X,0)}$ is the natural lift of $X$ from $M$ to
$TM\oplus T^*M$, while $Z_{(0,\xi)}$, at a point $(Y,\upsilon)\in
TM\oplus T^*M$, is equal to the vertical vector with value $-i_Y d\xi$.

If $\pi$ is a bivector field on $M$, the graph $L_{\pi}$ of
$\tilde{\pi}:T^*M\to TM$ is an $E_0$-Dirac structure if and only if
$\pi$ is a Poisson structure, i.e. if  $[\pi,\pi]=0.$
The Lie algebroid structure on $L_{\pi}$ may be transferred by
projection to $T^*M$, where it becomes the usual Lie algebroid
structure associated to $\pi$.

For a 2-form $\omega$ on $M$, the graph $L_{\omega}$ of
$\tilde{\omega}:TM\to T^*M$ is an $E_0$ Dirac structure if and only if
$d\omega=0$.  When the Lie
algebroid structure on $L_{\omega}$ is transferred to $TM$ by
projection, it always becomes the standard one.

\section{Twisted Poisson structures}

Now let $\phi$ be a 3-form on $M$.  We define a new
bracket on $E_0$ by adding the term $\phi(X_1,X_2,\cdot)$
to the right hand side of (\ref{eq-courant}).  A
simple computation (see the beginning of
Section \ref{sec-gauge} for another argument using less computation) shows
that the new bracket together
with the original bilinear form and anchor constitute a Courant
algebroid structure on $TM\oplus T^*M$ if and only if $d\phi=0$.
We denote this modified Courant algebroid by $E_{\phi}$.

We call a 2-form $\omega$ on $M$ {\bf $\phi$-closed} if
$L_{\omega}$ is an $E_{\phi}$ Dirac structure.  This just means
that $d\omega=\phi$; hence, $\phi$-closed forms exist only when
$\phi$ is exact.  Of course, any 2-form becomes
``twisted closed'' if we allow ourselves to choose $\phi$ ``after the
fact.''

More interesting are the {\bf $\phi$-Poisson structures},
i.e. bivector fields $\pi$ for which $L_{\pi}$ is an
$E_{\phi}$-Dirac structure.  It is easily seen that $\pi$ is a
$\phi$-Poisson structure if and only if it satisfies (\ref{eq-quasi}).

To understand the meaning of (\ref{eq-quasi}), we introduce
Poisson brackets and hamiltonian vector fields by the usual definitions; i.e.
$\{f,g\}=\pi(df,dg)$ and $H_f=\{\cdot,f\}.$  The usual Jacobi
equation then acquires an extra term:
\begin{equation}
\label{eq-jacobi}
\{\{f,g\},h\} + c.p.  + \phi(H_f,H_g,H_h) = 0,
\end{equation}
where ``$c.p.$'' means the sum of the two terms obtained from the
previous expression by circular permutation of the three variables.

If $M$ is 3-dimensional, the hamiltonian vector fields span a space of
dimension 0 or 2 at each point, so the extra term in the Jacobi identity is
zero; i.e. a twisted Poisson structure is a Poisson structure.
Nevertheless, as we shall see, the presence of $\phi$ still has an
effect on some of the standard Poisson-geometric constructions.  Thus,
one should consider $\phi$ as part of the twisted Poisson structure,
even if it does not affect the Jacobi identity.

Twisted Poisson structures
lead to Lie algebroid structures on the cotangent bundle,
just as ordinary Poisson structures do, since any
Dirac structure is a Lie algebroid for the induced
Courant bracket; identifying $L_{\pi}$ with $T^*M$ by
projection to the second summand of $E_{\phi}$ transfers the Lie
algebroid structure to the cotangent bundle.  The explicit formula for this
bracket is:
\begin{equation}
\label{eq-formbracket}
[\omega_1,\omega_2] = \call_{\tilde\pi(\omega_1)}\omega_2 -
\call_{\tilde\pi(\omega_2)}\omega_1 - d(\pi(\omega_1,\omega_2))
+ \phi(\tilde{\pi}(\omega_1),\tilde{\pi}(\omega_2),\cdot).
\end{equation}
Since this bracket satisfies the Jacobi identity, but the Poisson
bracket of functions does not, it should not be surprising that the relation
between the two brackets is also altered by the presence of
$\phi$, namely
\begin{equation}
\label{eq-twobrackets}
[df,dg]=d\{f,g\}+\phi(H_f,H_g,\cdot).
\end{equation}

Similarly, the mapping from functions to their hamiltonian vector
fields is no longer an antihomomorphism.  Instead, we have
$$H_{\{f,g\}}+[H_f,H_g]=-\tilde\pi(\phi(H_f,H_g,\cdot)).$$

The dual of the Lie algebroid $L_\pi$ may be identified with the
tangent bundle $TM$.  The Lie algebroid cohomology differential is no
longer simply $d_\pi = [\pi,~]$
(which is not of square zero unless
$[\pi,\pi]=0$), but is now given by the formula: $$d_{\pi,\phi}=d_\pi+(\wedge^2\tilde\pi
\otimes 1)(\phi) .$$  In this formula, $(\wedge^2\tilde\pi
\otimes 1)(\phi)$ is a section of the tensor bundle $\wedge^2 TM
\otimes T^*M$.  It acts as a degree 1 operator on multivector fields
by contraction with the factor in $T^*M$.  This operator is zero
on functions, so for $f\in C^\infty(M)$ we have
$$d_{\pi,\phi}f=d_\pi f = [\pi,f]=H_f$$ as usual.  On the other
hand, for a vector field $X$,
$$d_{\pi,\phi}X=d_\pi X +(\wedge^2\tilde\pi
\otimes 1)(\phi) X,$$ which operates on a pair of 1-forms by
$$(d_{\pi,\phi}X)(\omega_1,\omega_2)=-(\call_X\pi)(\omega_1,\omega_2)+
\phi(\tilde\pi\omega_1,\tilde\pi\omega_2,X).$$
Since $d_{\pi,\phi}^2=0,$ we have for each hamiltonian vector field
$H_f$
$$(\call_{H_f}\pi)(\omega_1,\omega_2)=\phi(\tilde\pi\omega_1,\tilde\pi\omega_2,H_f),$$
so the flow of a hamiltonian vector field does not in general
preserve $\pi$.  This is, of course, evident from the failure of
the Jacobi identity.  Incidentally, hamiltonian flows do not in
general preserve $\phi$ either, as may be seen easily in the
3-dimensional case, where $\phi$ is arbitrary.

An important feature of Poisson manifolds is their decomposition
into symplectic leaves.  For twisted Poisson manifolds, we have a
similar decomposition, into the orbits of the Lie algebroid
$L_\pi$ (whose anchor is $\tilde\pi$, just as in the ordinary
case).  Each orbit carries a nondegenerate 2-form, but these forms
are only twisted closed; their differentials equal the pullback of
the 3-form; i.e. a twisted Poisson manifold is decomposed into
twisted symplectic leaves.

\section{Gauge transformations associated to 2-forms}
\label{sec-gauge}
For any 2-form $B$ on $M$, we define the endomorphism
$\tau_B$ of $TM\oplus T^*M$ by
$\tau_B(X,\xi)=(X,\xi+\tilde{B}(X))$.  A simple computation shows that
$\tau_B$ preserves the symmetric bilinear form and the anchor common
to all the $E_{\phi}$  and that, for any closed 3-form $\phi$, $\tau_B$ is an isomorphism
of Courant bracket structures from $E_\phi$ to $E_{\phi-dB}$.
This shows that $\tau_B$ is an automorphism of $E_\phi$ if and only if
$B$ is closed, and that
the isomorphism class of $E_\phi$
depends only on the class $[\phi]$ in the de Rham cohomology space
$H^3(M,\reals)$.  Since any closed 3-form $\phi$ is locally exact,
application of a transformation $\tau_B$
 also shows immediately that the bracket for $E_\phi$
does indeed satisfy the Courant algebroid
axioms.

If $\xi$ is any 1-form, then the gauge transformation $\tau_{-d\xi}$ is
the time-1 flow of $Z_{(0,\xi)}$.  This recovers without further
computation the fact that $\tau_{-d\xi}$ is an automorphism; since any
closed 2-form is locally exact, it also gives another proof of the
fact that gauge transformation by a closed 2-form is an automorphism.

The additive group of closed 2-forms on $M$ acts on the space
of $\phi$-Dirac structures for each $\phi$, while
if $\phi=dB$, $\tau_{B}$ maps $\phi$-Dirac structures to ordinary
Dirac structures.  We will refer to these operations (whether $B$ is
closed or not) as {\bf gauge
transformations}.  Furthermore, when gauge-equivalent Dirac
structures correspond to (twisted) closed 2-forms or (twisted) Poisson
structures, we will say that the forms or Poisson structures are
gauge-equivalent.

Although gauge-equivalence looks like a rather coarse relation (for
instance, any two closed 2-forms are gauge-equivalent as $E_0$-Dirac
structures), it is still nontrivial.  First of all, the Lie algebroid
structures on gauge-equivalent Dirac structures are always isomorphic,
since the Lie algebroid operations on any $E$-Dirac structure are
simply the restriction of the Courant algebroid bracket and anchor on
$E$.

Gauge-equivalence is especially interesting when applied to
twisted Poisson structures.  Given a bivector $\pi$ and a 2-form $B$,
$\tau_B(L_\pi)$ is of course a maximal isotropic subspace of
$TM\oplus T^*M$, but it is not necessarily of the form
$L_{\pi'}$ for another bivector $\pi'$.  If it is, then
$\tilde{\pi'}=\tilde\pi(1+\tilde{B}\tilde\pi)^{-1}.$
(This formula is equivalent to $\tilde{\pi'}^{-1}=\tilde\pi^{-1}+
\tilde B$ when  $\tilde{\pi'}$ and $\tilde\pi$ are
invertible.)
In fact, the
condition for $\tau_B(L_\pi)$ to correspond to a bivector is
just that the endomorphism $1+\tilde{B}\tilde\pi$ of $TM$ be
invertible.  By abuse of notation, we will write   $\tau_B\pi$ for
$\pi'$.  For general $\pi$ and $B$, we may think of $\tau_B\pi$
as a bivector with singularities.

Since gauge-equivalent twisted Poisson structures $\pi$ and $\tau_B \pi$
give rise to isomorphic
Lie algebroids, they have many features in common.  (Note that
 the particular Lie algebroid
structure on $T^*M$ is changed under a gauge transformation of
$\pi$, though the structures coming from $\pi$ and $\tau_B \pi$
are isomorphic via the map $1+\tilde B \tilde \pi$.)
Since the (twisted) Poisson cohomology for $\pi$ is isomorphic to the
Lie algebroid cohomology of $L_\pi$, (twisted) Poisson cohomology is
gauge invariant.
Gauge-equivalent structures have the same Casimir
functions.  Their
decompositions into twisted symplectic leaves are the same, though
the 2-forms along the leaves differ by the pullbacks of $B$.  In
particular, for $\phi=0$ and $dB=0$, the variation from leaf to
leaf of the cohomology class of the symplectic structure along the
leaves (the fundamental class of Dazord and Delzant
\cite{da-de:probleme}) is the same.

It is natural to ask at this point whether every twisted Poisson
structure is gauge-equivalent to an ordinary Poisson structure.
If $\pi$ is a $\phi$-Poisson structure and $dB=\phi$, then
$\tau_{-B}\pi$ will exist and be a Poisson structure if and only
if $1-\tilde B\tilde\pi$ is invertible.  (If it is not,
$\tau_{-B}\pi$ still makes sense as a Dirac structure; see
Examples
\ref{ex-global} and \ref{ex-group} below.)
So the problem is to find a primitive $B$ for
$\phi$ such that $1-\tilde B\tilde \pi$ is invertible.
If $\phi$ is not exact, we have no hope of finding a global
primitive $B$; Example \ref{ex-global} shows even an exact $\phi$ cannot
always be ``gauged away.''

Using the Poincar\'e lemma, it is always possible to find a locally
finite covering of $M$ by open subsets $\calu_{\alpha}$ carrying
primitives $B_\alpha$ for $\phi$ which are small enough so that
$1-\wt{ B_\alpha} \tilde \pi$ is invertible.  We thus have a family
of locally defined ordinary Poisson  structures
$\pi_\alpha=\tau_{-B_\alpha}\pi$ which are related on overlaps
$\calu_\alpha\cap \calu_\beta$ by gauge transformations
corresponding to the exact \v{C}ech 1-cocycle $B_\alpha - B_\beta$
with values in 2-forms. The usual argument due to Weil
\cite{we:sur} takes us from here to a real-valued 3-cocycle
representing the cohomology class of $\phi$.  In this way, our
twisted Poisson structures resemble other twisted structures, such
as vector bundles, arising from the presence of cohomology classes
of degree 3.  Although  we are now veering very close to the world
of gerbes which we have tried to avoid, we will return to this
point in the last section.

\section{Examples}

The following example shows that, even when $\phi$ is exact, it
might not be possible to remove it by a global gauge
transformation of twisted Poisson structures.  Since our example will
be 3-dimensional, we remind the reader that $\phi$ should be
considered as part of the structure, even if $\pi$ already satisfies
the Jacobi identity.

\vskip 3mm
\begin{ex}
\label{ex-global}
{\em
On $\reals^3$, we consider the Lie-Poisson structure
$\pi=x_1\frac{\del}{\del x_2}\wedge\frac{\del}{\del x_3} + c.p. $
as a quasi-Poisson structure with 3-form $\phi=3\lambda  dx_1\wedge dx_2\wedge
dx_3,$ where $\lambda $ is an arbitrary real constant.  An obvious primitive
for $\phi$ is $B=\lambda (x_1 dx_2\wedge dx_3 + c.p.);$  when we apply the
corresponding gauge transformation, we find that the Poisson structure
$\tau_{-B}\pi$ is equal to $\pi/(1+\lambda r^2),$ where
$r^2=x_1^2+x_2^2+x_3^2$.
  If $\lambda$ is
positive, this is fine, but if $\lambda$ is negative, the structure is
singular along the sphere of radius $(-\lambda)^{-1/2}$.  The
corresponding Dirac structure has this sphere as a ``presymplectic
leaf'' carrying the zero 2-form.

We may ask whether a gauge transformation using
 another choice of primitive would lead to a
globally defined Poisson structure, i.e. whether we can have
$1-\tilde B\tilde\pi$ invertible.  In fact, this is not possible.  To
 see this, we use the usual identifications arising from the orientation and
 metric on $\reals^3$ to replace the skew-symmetric $3\times 3$
 matrix-valued
 functions $\tilde
 B$ and $\tilde\pi$ by 3-dimensional vector fields $\mathbf B$ and
$\mathbf x$.  (The latter is just the identity
 vector field.) The conditions which must be satisfied by
 $\mathbf B$ are ${\mathbf\nabla}\cdot{\mathbf B} = 3 \lambda$ and ${\mathbf x}
 \cdot {\mathbf B} ({\mathbf x}) \neq -1$.  Since $\mathbf x$ vanishes at the
 origin, we must have  ${\mathbf x} \cdot {\mathbf B} ({\mathbf x}) > -1$.  Taking
 the surface integral of $\mathbf B$ over a sphere of radius $r$ and
 applying the divergence  theorem, we obtain the inequality $\lambda
 r^2 > -1$.  If $\lambda$ is negative, this means that $\tau_{-B}\pi$
 cannot be a nonsingular Poisson structure beyond the sphere of radius
$(-\lambda)^{-1/2}$.
}
\end{ex}
\vskip 3mm

The next example was first found (in a more general form) by
\Klimcik ~and \Severa
 \cite{kl-se:wznw}; it was rediscovered independently by Alekseev
and Str\"obl (unpublished).

\vskip 3mm
\begin{ex}
\label{ex-group}
{\em
On a Lie group $G$ with bi-invariant metric (not necessarily definite), a natural
choice of 3-form is the bi-invariant Cartan form (sometimes called the
Chern-Simons form) defined on the Lie algebra by
$\phi(u,v,w)={\textstyle{\frac 12}}[u,v]\cdot w$. We will use the metric to identify
$TG$ with $T^*G$; when $\phi$ is thus interpreted as a bilinear map
on vector fields with values in vector fields, its value on a pair
of left-invariant and right-invariant vector fields is just half their
bracket.  For any element $a$ of the Lie algebra $\frakg$, we will
denote the corresponding left-invariant and right-invariant vector
fields by $a^L$ and $a^R$.

In $E=TG\oplus TG$, we consider the maximal isotropic subbundle $L$
given by the values
of the sections $e_a=(a^L-a^R,{\textstyle{\frac 12}}(a^L +a^R))$ as $a$
ranges over $\frakg$.  Computing the ordinary Courant bracket of
$e_a$ and $e_b$ gives an expression which differs from $e_{[a,b]}$
by a term in the 1-form component which vanishes if we add
$-\phi(e_a,e_b,\cdot)$ to the bracket; i.e.
$L$ is an $E_{-\phi}$-Dirac structure.

Since the map $a\mapsto e_a$ from $\frakg$ to sections of $L$
 is an isomorphism to each fibre of $L$, it follows that the
Lie algebroid $L$ is isomorphic to the action Lie algebroid.
Examining the component of $e_a$ in $TG$, we see that the action
is just the conjugation action, so the orbits of the Lie
algebroid, and hence the twisted presymplectic leaves
of the twisted Dirac structure, are the connected components of
the conjugacy classes in $G$.

$L$ is a twisted Poisson structure over the subset of $G$
where the vector fields $a^L+a^R$ are linearly independent, i.e.
the open, dense subset
where $\Ad_g+1$ is invertible.  Over this subset, $L$ is of the
form $L_\pi$, where $\tilde\pi(g)$ is the conjugate via left
translation by $g$ of the operator $2(\Ad_g-1)/(\Ad_g+1)$  on
$\frakg$.  (Recall that we are identifying tangent and cotangent
vectors with the inner product.)  As in the previous example, we
may also wish to  think of $\pi$ as a singular twisted Poisson structure on the
whole group; on the conjugacy classes where it blows up, the
twisted closed 2-form becomes degenerate.

This example also shows the utility of computing with Dirac
structures.  One could define $\pi$ directly in terms of the operator
$2(\Ad_g-1)/(\Ad_g+1)$, but it would be much more complicated to
compute $[\pi,\pi]$ than it is to compute the Courant bracket of
elements of $L_\pi$.

We end this example by noting that, if one replaces the operator $ 2(\Ad_g-1)/(\Ad_g+1)$
by $\Ad_g$ itself, one obtains a globally defined bivector field $\pi_1$ which is a
{\em quasi}-Poisson structure in the  sense of Alekseev {\it et al.}~
\cite{al-ko-me:quasi};
i.e. $[\pi_1,\pi_1]$ is a multiple of the Cartan trivector field
associated to the Cartan 3-form by the invariant metric.  Its ``quasi-symplectic
leaves'' are once again the connected components of the conjugacy classes.
It would be very interesting to understand the relation between the two
``nearly Poisson'' structures on $G$.
}
\end{ex}

\section{Deformation quantization and symplectic groupoids}

This final section is devoted to some speculations.

Poisson structures on $M$ correspond to infinitesimal deformations
of the associative algebra structure on $C^\infty (M)$; see Bayen
{\it et al.}~\cite{bffls:deformation} and Kontsevich
\cite{ko:deformation}.  To what algebraic objects do twisted
Poisson structures correspond?  As already suggested by Park \cite{pa:topological}
as well as by Cornalba
and Schiappa \cite{co-sc:nonassociative}, the anomaly in the Jacobi
identity for twisted Poisson structures should correspond to an
anomaly in associativity for deformations of $C^\infty (M)$, but
the precise algebraic nature of the nonassociativity has yet to be
determined. On the other hand, our approach to twisted Poisson
structures via gauge transformations suggests a different approach
to quantization.

As in the ordinary deformation quantization, we start with a formal $\phi$-Poisson
structure $\pi\in\Gamma(\bigwedge^2 TM)[[\epsilon]]$ such that $\pi=O(\epsilon)$,
i.e. $$\pi=\pi^{(1)}\epsilon +\pi^{(2)}\epsilon^2+\pi^{(3)}\epsilon^3+\cdots.$$
We choose  a
cover of $M$ by open subsets $\calu_\alpha$ carrying 2-forms
$B_\alpha$ such that $dB_\alpha=\phi$, so that $\tau_{-B_\alpha}\pi=\pi_\alpha$ become
(formal) Poisson structures.  Assuming
further that we have chosen a quasi-isomorphism of Kontsevich for $M$, we obtain
an associative $\reals[[\epsilon]]$-algebra
structure on $\cala_\alpha=C^\infty(\calu_\alpha)[[\epsilon]]$.

On $\calu_\alpha\cap \calu_\beta$, we have $\pi_\alpha=\tau_{
B_\alpha-B_\beta}\pi_\beta.$ We may write $B_\alpha-B_\beta =
d\theta_{\alpha\beta}$, where the $\theta_{\alpha\beta}$ form a
1-cochain of 1-forms which is a cocycle only when considered modulo
closed 1-forms.  The exactness of $B_\alpha-B_\beta$ implies that
$\pi_\alpha$ and $\pi_\beta$ are equivalent Poisson structures,
i.e. there is a formal curve of diffeomorphisms carrying
$\pi_\alpha$ to $\pi_\beta$. To prove this we construct a
time-dependent formal Poisson structure $\pi(t)$ connecting
$\pi_\alpha$ with $\pi_\beta$, namely
$\pi(t)=\tau_{t(-B_\alpha+B_\beta)}\pi_\alpha$, and a formal
time-dependent vector field $v(t)$, $v(t)=O(\epsilon)$, such that
$d\pi(t)/dt={\cal L}_{v(t)}\pi(t)$. It is easy: we know that
$\tau_{t(-B_\alpha+B_\beta)}$ is the flow of
$Z_{(0,\theta_{\alpha\beta})}$, but we know that the
flow of $Z_{(\wt{\pi(t)}\theta_{\alpha\beta},\theta_{\alpha\beta})}$
preserves $L_{\pi(t)}$, so that we can put
$v(t)=-\wt{\pi(t)}\theta_{\alpha\beta}$.  By formality there are
isomorphisms $I_{\alpha\beta}$ between the quantized algebras on the
overlap regions.  On a triple intersection $\calu_\alpha \cap
\calu_\beta \cap \calu_\gamma$, the product
$I_{\alpha\beta}I_{\beta\gamma}I_{\gamma\alpha}$ cannot be expected to
be the identity; rather it should be an inner automorphism related to
a function $f_{\alpha\beta\gamma}$ which is a primitive for the closed
1-form
$\theta_{\alpha\beta}+\theta_{\beta\gamma}+\theta_{\gamma\alpha}.$ On
quadruple intersections, these functions in turn will give rise to a
\v{C}ech 3-cochain which represents the cohomology class $[\phi]$.

The exact nature of the resulting algebraic structure and its relation
to nonassociativity are still unclear to us.

Let us remark here (related constructions were considered by Park \cite{pa:topological})
that the graded Lie algebra of multivector fields
$\Gamma(\bigwedge TM)[1]$ is twisted by $\phi$ to an $L_\infty$
algebra denoted $\Gamma(\bigwedge TM)[1]_\phi$, such that formal
$\phi$-Poisson structures are the formal solutions of Maurer-Cartan
equation.  The transformations $\tau_B$ give rise to
$L_\infty$-isomorphisms $\Gamma(\bigwedge
TM)[1]_\phi\rightarrow\Gamma(\bigwedge TM)[1]_{\phi-dB}$, hence
$\Gamma(\bigwedge TM)[1]_\phi$ can be seen as glued from the graded
Lie algebras $\Gamma(\bigwedge T{\cal U}_\alpha)[1]$. It is natural to
ask if there is a natural twist of the Hochschild complex that would
be quasi-isomorphic to $\Gamma(\bigwedge TM)[1]_\phi$ (of course, we
can define a twist just by using Kontsevich's quasi-isomorphism).

To describe the $L_\infty$-structure on $\Gamma(\bigwedge
TM)[1]_\phi$, we can use the language (and ideas) of Kontsevich 
\cite{ko:deformation}.  We define a graded $Q$-manifold $X_\phi$:
as a graded supermanifold, $X_\phi =T[1]M\times\reals[2]$, the
vector field $Q$ is $d+\phi\partial_t$, where $d$ is the deRham
differential on $M$ and $t$ is the coordinate on $\reals[2]$. Let
$j^1X_\phi$ be the space of one-jets of sections of the projection
$X_\phi\rightarrow T[1]M$; it is again a graded $Q$-manifold and so
is the (infinite-dimensional) space of legendrian submanifolds of
$j^1 X_\phi$, denoted $\Lambda_\phi$. Then the formal graded
$Q$-manifold corresponding to $\Gamma(\bigwedge TM)[1]_\phi$ is
just the formal neighborhood of a point in $\Lambda_\phi$ --- the
legendrian submanifold of $j^1X_{\phi}$ given by the partially defined
section of $X_\phi\rightarrow T[1]M$, $m\mapsto (m,0)$, $m\in M$.
The transformations $\tau_B$ act on $X_\phi$ by $(y,t)\mapsto
(y,t+B)$; this action is inherited by $\Lambda_\phi$ and
$\Gamma(\bigwedge TM)[1]_\phi$.

Finally, we recall that symplectic groupoids provide a geometric
model (Weinstein \cite{we:noncommutative}) as well as a potential
means of construction (Cattaneo and Felder \cite{ca-fe:poisson})
for deformed algebras.  The twisted version of this story appears
to go roughly as follows, following the pattern of
Kosmann-Schwarzbach \cite{ko:jacobian}, which treated the case
where $M$ is a point.  Given a $\phi$-Poisson structure $\pi$ on
$M$, the $\phi$-Dirac structure $L_\pi$ and the ``horizontal''
subbundle $TM$ form a Manin pair in $E_\phi$, and so they form a
``quasi-Lie bialgebroid.''  Assuming that the Lie algebroid
$L_\pi$ can be integrated to a groupoid $G$, an extension of the
methods of Mackenzie and Xu \cite{ma-xu:integration} should permit
the construction of a nondegenerate twisted Poisson structure on
$G$ which is compatible with the groupoid structure, making $G$
into a ``twisted symplectic groupoid.'' (The twisting 3-form would
be the difference between the pullbacks of $\phi$ from $X$ to $G$
by the target and source maps.) The relation of this kind of object
to deformations remains a mystery.

\end{document}